\documentclass[11pt]{article}
\usepackage{amsmath}
\usepackage{amssymb}
\usepackage{latexsym}
\usepackage{graphics}
\usepackage[all]{xy}

\usepackage{url} 

\usepackage{graphicx}
\usepackage{float}

\usepackage{amsthm}

\newtheorem{df}{Definition}[section]
\newtheorem{thm}[df]{Theorem}
\newtheorem{prop}[df]{Proposition}

\newtheorem{lem}[df]{Lemma}

\newcommand{\pfend}%
{\hspace*{\fill}\lower3pt\hbox{$\Box$}\medskip}

\newcommand{\R}{\mathbb{R}}
\newcommand{\Z}{\mathbb{Z}}
\newcommand{\C}{\mathbb{C}}
\newcommand{\qt}{\mathbb{H}}
\newcommand{\RP}{\mathbb{RP}}
\newcommand{\CP}{\mathbb{CP}}
\newcommand{\HP}{\mathbb{HP}}

\begin{document}

\title{\textbf{The Gromov-Lawson-Rosenberg conjecture for the Semi-dihedral group of order 16}}
\author{Arjun Malhotra and Kijti Rodtes\\
        }

\maketitle
\section{Introduction}
The question of which manifolds admit metrics of positive scalar curvature is an extremely interesting and open one. By a well known result of Lichnerowicz \cite{12}, there are spin manifolds which do not admit such metrics. Indeed, by the Lichnerowicz formula, the existence of such a metric implies that the index of the Dirac operator vanishes. This, combined with the Atiyah-Singer index theorem implies that a topological invariant known as the $\hat{A}$ genus, which is a linear combination of the Pontrjagin classes of the manifold, vanishes.

The $\hat{A}$ obstruction was generalized by Hitchin \cite{7} to an obstruction $\alpha(M) \in KO_{*}$. Letting $\pi$ denote any fundamental group, this gives rise to a transformation of cohomology theories $\alpha:\Omega ^{spin}_{*}(B\pi) \to KO_{*}(B\pi)$, and Gromov and Lawson conjectured that $\alpha(M)=0$ was also a sufficient condition for $M$ to admit a metric of positive scalar curvature.

Rosenberg \cite{ro3} later generalized this further to define an obstruction $ind:\Omega ^{spin}_{*}(B\pi) \to KO_{*}(C^{*}_{red}\pi)$, thought of as an equivariant generalized index.  Modifying the Gromov-Lawson conjecture, Rosenberg conjectured that the converse was true also; namely that a compact spin manifold $M^{n}$ with $\pi_{1}(M)=\pi$ and $n \geq 5$ admits a positive scalar curvature metric if and only if $ind[u:M\to B\pi]=0 \in KO_{n}(C^{*}_{red}\pi)$ where $ind:\Omega ^{spin}_{*}(B\pi) \to KO_{*}(C^{*}_{red}\pi)$ can be thought of as an equivariant generalized index map, and $u$ is the classifying map for the universal cover of $M$.

The conjecture has been proven in the simply connected case \cite{17}, if $\pi$ has periodic cohomology \cite{2}, and if $\pi$ is a free group, free abelian group, or the fundamental group of an orientable surface \cite{15}. It is now known to be false in general, for example if $\pi=\Z^{4} \times \Z_{3}$, and for a large class of torsion free groups \cite{16}, \cite{4}.

The aim of this paper is to prove the following result:
\begin{thm}
The Gromov-lawson-Rosenberg conjecture (GLR) is true for the semi-dihedral group of order 16.
\end{thm}
We recall that the index map factors through connective and periodic real K-theory as follows:
$$ind=A \circ p \circ D: \Omega_{n}^{spin}(B \pi) \rightarrow ko_{n}(B \pi)\rightarrow KO_{n}(B \pi) \rightarrow KO_{*}(C^{*}_{red}\pi). $$
Here $D$ is a generalized Dirac operator, which sends a spin manifold to its $ko$ fundamental class, $p$ is the periodicity map, and $A$ is known as the assembly map.  We denote by $ko_{n}^{+}(B \pi)$ the set $D(\Omega_{n}^{spin,+}(B \pi))$, where $ \Omega_{n}^{spin,+}(B \pi)$ is the subgroup represented by pairs $(N,f)$ for which $N$ admits a positive scalar curvature metric. The following result, due to Jung and Stolz, is the basis of our proof.
\begin{thm}
A compact spin manifold $M^{n}$ with $\pi_{1}(M)=\pi$ and $n \geq 5$ admits a positive scalar curvature metric if and only if $D(M,u) \in ko_{n}^{+}(B \pi)$, where $u$ is the classifying map for the universal cover of $M$.
\end{thm}
Thus one way of proving the conjecture is to first calculate $ko_{n}(B \pi)$, and then to identify $Ker(Ap)$ and realize all of it by positive scalar curvature manifolds, and this is what we shall do in the case of the semi-dihedral group $SD_{16}$

The next section describes the calculations of $ko_{*}(BSD_{16})$. The method used is the Greenlees spectral sequence.

The following two sections then focus on realizing $Ker(Ap) \subset ko_{*}(BSD_{16})$. This is done by detection in periodic K-theory and in ordinary homology. Firstly, we consider inclusions of lens spaces from cyclic and quaternion subgroups, and follow the method of calculating eta-invariants in \cite{2}. We show that this detects all of $Ker (A) \cap Im(p)$.

This leaves certain Bott torsion classes to detect, and this is done by first considering inclusions of projective bundles over projective space from a Klein-$4$ subgroup. The method here is to detect these classes in ordinary homology, by making explicit calculations. The result is that the classes appearing in second local cohomology may all be realized by this inclusion.

This leaves the Bott torsion in first local cohomology. The eta-multiples are immediately dealt with, since the subgroup spanned by fundamental classes of positive scalar curvature manifolds is an ideal, but it turns out that there is a class appearing in dimensions $8k$ which comes from local cohomology of the eta multiples in $ko^*(BSD_{2^N})$, but which is nevertheless not an eta-multiple itself. Calculation shows that this class is embedded in the ordinary homology $H_{8k}(BSD_{2^N}; \Z_2)$, with fundamental class \emph{distinct} from those realized by the inclusion from the Klein-$4$ subgroup above. We realize this class by giving an explicit geometric construction of a lens space bundle over a circle (compare \cite{8}).

\section{The Greenlees spectral sequence and ko calculations}
In this section we show how to calculate $ko_*(BSD_{16})$. For a group $G$, we can calculate $ko_{*}(BG)$ in several ways. For instance, by the Adams spectral sequence with input $H_{*}(BG;\mathbb{F}_{2})$, by the Atiyah-Hirzebruch spectral sequence or by \emph{Bruner-Greenlees method}.  Here, \emph{Bruner-Greenlees method} can be described by the diagram below;\\
\setlength{\unitlength}{0.9cm}
\begin{picture}(10,4)
\linethickness{0.2mm}

\put(1,1.75){$H^{*}(BG;\mathbb{F}_{2})$}
\put(5,3){$ku^{*}(BG)$}
\put(10,3){$ku_{*}(BG)$}
\put(5,0.5){$ko^{*}(BG)$}
\put(10,0.5){$ko_{*}(BG),$}

\put(3.1,2){\vector(2,1){1.8}}
\put(6.7,3.1){\vector(1,0){3.2}}
\put(6.7,0.6){\vector(1,0){3.2}}
\put(5.7,2.8){\vector(0,-1){1.8}}
\put(10.7,2.8){\vector(0,-1){1.8}}
\put(3.3,2.5){ASS}
\put(8,3.2){GSS}
\put(8,0.7){GSS}
\put(5.8,1.8){BSS}
\put(10.8,1.8){BSS}
\end{picture}\\
where ASS refer to the Adams spectral sequence, BSS refer to the $\eta$-Bockstein spectral sequence and GSS refer to the Greenlees spectral sequence.  From this diagram, we see that to obtain $ko_{*}(BG)$ via Bruner-Greenlees method with the input $H^{*}(BG;\mathbb{F}_{2})$ we can proceed in two ways around the square, namely $ASS\longrightarrow BSS\longrightarrow GSS$ and $ASS \longrightarrow GSS \longrightarrow BSS$.  The first way is suitable for tackling the Gromov-Lawson-Rosenberg conjecture, since we have:
\begin{lem}\label{lem Ap}(\cite{3}, Lemma 2.7.1.) The image of $Ap$ is isomorphic to the 0-column at $E_{\infty}$-page of the Greenlees spectral sequence for $ko_{*}(BG)$.
The kernel of $Ap$ has a filtration with subquotients given by the higher columns
at $E_{\infty}$-page.
\end{lem}

The calculation of $ku^{*}(BG)$ for a finite group $G$ from  $H^{*}(BG;\mathbb{F}_{2})$ may be performed using the Adams spectral sequence:
 $$E^{s,t}_{2}=\operatorname{Ext}^{s,t}_{E(1)}(\mathbb{F}_{2},H^{*}(BG;\mathbb{F}_{2}))\Longrightarrow ku^{-(t-s)}(BG)_{2}^{\wedge},$$
where $E(1)$ denotes the exterior algebra on the Milnor generators $Q_{0}$ and $Q_{1}$, and more details may be found in \cite{BGku}.  The calculation of $ko^{*}(BG)$ from $ku^{*}(BG)$ by using BSS;
 $$ E^{*,*}_{1}=ku^{*}(BG)[\widetilde{\eta}] \Rightarrow ko^{*}(BG),$$
where $\widetilde{\eta}$ has bidegree $(1,1)$ and the differential $d_{r}:E^{s,t}_{r}\longrightarrow E^{s+r,t-1}_{r}$, is recently developed by R.R. Bruner and J.P.C. Greenlees in \cite{3}.  To gain $ko_{*}(BG)$, as stated above, we prefer to use the Greenlees spectral sequence, GSS \cite{3};
$$E^{s,t}_{2}=H^{-s}_{I}(ko^{*}(BG))_{t}\Rightarrow  ko_{(s+t)}(BG),$$ where $I=\ker(ko^{*}(BG)\longrightarrow ko^{*})$, the augmentation ideal of $ko^{*}(BG)$, and differentials $d_{r}:E^{s,t}_{r}\longrightarrow E^{s-r,t+r-1}_{r}$.  

The strategy of GSS for $ko_{*}(BG)$ is mainly the decomposition of the input, $ko^{*}(BG)$, via the short exact sequences;
\begin{equation}\label{eqre1}
    0\longrightarrow ST \longrightarrow \xymatrix{ko^{*}(BG)\ar[r]^{\pi^{o}}& QO}\longrightarrow0
\end{equation}
and
\begin{equation}\label{eqre2}
    0\longrightarrow T \longrightarrow \xymatrix{ko^{*}(BG)\ar[r]^{\pi^{u}} &\overline{QO}}\longrightarrow0
\end{equation}
where $QO$ is the image of $ko^{*}(BG)$ in $KO^{*}(BG)$ and $\overline{QO}$ is the image of $ko^{*}(BG)$ in $KU^{*}(BG)$ such that $ST$ is the $\beta$-torsion part of $ko^{*}(BG)$ and $T$ is the $\ker \pi^{u}$.  Here, $QO$ and $\overline{QO}$ are modules over $ko^{*}(BG)$ via $\pi^{o}$ and $\pi^{u}$ respectively.
Moreover, let $Q\tau$ be the kernel $\ker(i:QO\longrightarrow \overline{QO})$.  By the snake lemma, $Q\tau\cong \operatorname{coker}(i':ST\longrightarrow T)$ or in other words, we have a short exact sequence
\begin{equation}\label{eqre3}
    0\longrightarrow ST \longrightarrow T\longrightarrow Q\tau\longrightarrow 0.
\end{equation}
Generally, $T$ will be $2$-torsion. 

Thus, to find $H^{*}_{I}(ko^{*}(BG))$, it is convenient to calculate from the long exact sequence induced by (\ref{eqre2}) together with the long exact sequence induced by (\ref{eqre3}). This is convenient because the calculation of $H_{I}^{*}(\overline{QO})$ involves character theory and that of $H_{I}^{*}(T)$ use mainly the ordinary cohomology ring.  More precisely, to get $H^{*}_{I}(ko^{*}(BG))$, we instead calculate $H_{I}^{*}(\overline{QO})$, $H_{I}^{*}(T)$ and the connecting homomorphism induced by (\ref{eqre2}), where $H_{I}^{*}(T)$ is obtained by $H_{I}^{*}(ST)$, $H_{I}^{*}(Q\tau)$ and the connecting homomorphism induced by (\ref{eqre3}).  Normally, $Q\tau \subseteq \tau$ $:=$ $\eta$-multiples, since only $\eta \in KO^{*}=\mathbb{Z}[\alpha,\eta,\beta,(\beta)^{-1}]/(\eta^{2},2\eta,\eta\alpha,\alpha^{2}-4\beta)$ is sent to $0 \in KU^{*}$.  Hence, if there is no $\eta$-multiple obtained from $TU$,  the $v$-torsion part of $ku^{*}(BG)$, then $Q\tau=\tau$.

Note that, for finite groups $G$, $ko^{*}(BG)$ is a Noetherian graded ring, \cite{3}, and thus local cohomology  is suitably calculated via the stable Koszul complex.  The first task for the calculation is to find radical generators for the ideals of each part $\overline{QO}, Q\tau, ST$, and we found in many cases that the number of generators required is less than or equal to the $p$-rank of $G$.  The next task is to find local cohomology followed by differentials of spectral sequence which can be done by using the connective property of $ko$ and the module structure of local cohomology.  The last task is to solve the extension problems which occur along the induced long exact sequence and at $E_{\infty}$-page.

For the case $G=SD_{16}$, we have (from Section 5.1.2 in \cite{kijti}) $Q\tau=\tau$ and $ST=TO$, where $TO$, by definition in \cite{3}, consists of Bockstein $\infty$-cycles in $ZTU$ lying on the zero line of the BSS for $ko^{*}(BSD_{16})$.  By Corollary 5.2.5 and Section 5.4 in \cite{kijti},
\begin{equation}\label{local cohomology of sd16}
    \begin{array}{ccl}
     H^{0}_{I}(ko^{*}(BSD_{16})) & = & H^{0}_{I}(\overline{QO})\oplus H^{0}_{I}(\tau), \\
     H^{1}_{I}(ko^{*}(BSD_{16})) & = & H^{1}_{I}(\tau)\oplus \ker(\delta:H^{1}_{I}(\overline{QO})\longrightarrow H^{2}_{I}(TO)), \\
     H^{2}_{I}(ko^{*}(BSD_{16})) & = & \operatorname{coker}(\delta:H^{1}_{I}(\overline{QO})\longrightarrow H^{2}_{I}(TO)),\\
    H^{i}_{I}(ko^{*}(BSD_{16})) &=& 0 \quad \hbox{for all $i\geq 3$,}
   \end{array}
\end{equation}
where $\delta$ is the connecting homomorphism in the long exact sequence induced by (\ref{eqre2}), and moreover $E_{2}$-page = $E_{\infty}$-page, see the details of the calculations in \cite{kijti}.  The $E_{\infty}$-page of GSS for $ko_{*}(BSD_{16})$ are shown in Figure 2.1 below ($R:=ko^{*}(BSD_{16})$).

\setlength{\unitlength}{1cm}
\begin{center}
\begin{picture}(20,17)
\multiput(0.2,2)(0,0.5){28}%
{\line(1,0){12.8}}
\put(1.5,2){\line(0,1){14}}
\put(3,2){\line(0,1){14}}
\put(11,2){\line(0,1){14}}
\put(13,2){\line(0,1){14}}
\multiput(11.9,4.6)(0,0.5){3}%
{$0$}
\multiput(11.9,8.6)(0,0.5){3}%
{$0$}
\multiput(11.9,12.6)(0,0.5){3}%
{$0$}

\multiput(11.9,3.6)(0,4){3}%
{$0$}
\multiput(11.8,2.6)(0,0.5){2}%
{$2^{3}$}
\multiput(11.8,6.6)(0,0.5){2}%
{$2^{4}$}
\multiput(11.8,10.6)(0,0.5){2}%
{$2^{5}$}
\multiput(11.8,14.6)(0,0.5){2}%
{$2^{5}$}
\multiput(11.8,2.1)(0,2){7}%
{$\mathbb{Z}$}

\multiput(6.8,2.6)(0,1){3}%
{$0$}
\multiput(6.8,7.6)(0,1){3}%
{$0$}
\multiput(6.8,11.6)(0,1){3}%
{$0$}
\multiput(6.8,3.1)(0,1){1}%
{$0$}

\put(6.8,2.1){$0$}
\put(6.8,5.6){$0$}
\put(5.6,4.1){$ [4]\oplus[8]\oplus[8]$}
\put(6.7,5.1){$[2]$}
\put(5.4,6.1){$[2]\oplus [4]\oplus[16]\oplus[32]$}
\put(6.7,6.6){$2$}
\put(6.3,7.1){$2\oplus [2]$}
\put(4.9,8.1){$[8]\oplus[16]\oplus[128]\oplus[128]$}
\put(6.7,9.1){$[4]$}
\put(4.7,10.1){$[2]\oplus [8]\oplus[16]\oplus[256]\oplus[512]$}
\put(6.6,10.6){$2^{2}$}
\put(6.2,11.1){$2^{2}\oplus [4]$}
\put(4.6,12.1){$[2]\oplus [32]\oplus[64]\oplus[2048]\oplus[2048]$}
\put(6.7,13.1){$[8]$}
\put(4.4,14.1){$[4]\oplus [32]\oplus[64]\oplus[16^{3}]\oplus[2\cdot16^{3}]$}
\put(6.6,14.6){$2^{2}$}
\put(6.2,15.1){$2^{2}\oplus [8]$}

\multiput(2,2.6)(0,1){13}%
{$0$}
\multiput(2,3.1)(0,1){1}%
{$0$}
\put(6.8,15.6){$\vdots$}
\put(2,15.6){$\vdots$}
\put(0.8,15.6){$\vdots$}
\put(11.8,15.6){$\vdots$}

\put(2,4.1){$0$}

\put(2,2.1){$0$}

\put(2,6.1){$0$}
\put(2,5.1){$2$}
\put(2,7.1){$2$}
\put(2,8.1){$2$}
\put(1.9,9.1){$2^{2}$}
\put(2,10.1){$2$}
\put(1.9,11.1){$2^{2}$}
\put(1.9,12.1){$2^{2}$}
\put(1.9,13.1){$2^{3}$}
\put(1.9,14.1){$2^{2}$}
\put(1.9,15.1){$2^{3}$}

\put(13.3,2.1){0}
\put(13.3,2.6){1}
\put(13.3,3.1){2}
\put(13.3,3.6){3}
\put(13.3,4.1){4}
\put(13.3,4.6){5}
\put(13.3,5.1){6}
\put(13.3,5.6){7}
\put(13.3,6.1){8}
\put(13.3,6.6){9}
\put(13.2,7.1){10}
\put(13.2,7.6){11}
\put(13.2,8.1){12}
\put(13.2,8.6){13}
\put(13.2,9.1){14}
\put(13.2,9.6){15}
\put(13.2,10.1){16}
\put(13.2,10.6){17}
\put(13.2,11.1){18}
\put(13.2,11.6){19}
\put(13.2,12.1){20}
\put(13.2,12.6){21}
\put(13.2,13.1){22}
\put(13.2,13.6){23}
\put(13.2,14.1){24}
\put(13.2,14.6){25}
\put(13.2,15.1){26}
\put(13.2,15.6){27}

\put(12.5,16.2){degree(t)}

\multiput(0.7,2.1)(0,0.5){27}%
{$0$}

\put(6.3,1.2){$H^{1}_{I}(R)$}
\put(11.5,1.2){$H^{0}_{I}(R)$}
\put(1.8,1.2){$H^{2}_{I}(R)$}
\put(0,1.2){$H^{\epsilon\geq3}_{I}(R)$}

\put(0.2,0.5){where$[n]:=$ cyclic group of order $n$, $2^{r}$:= elementary abelian group of rank $r$.}
\linethickness{0.5mm}
\put(11,2){\line(1,0){2}}
\put(11,2.5){\line(-1,0){8}}
\put(3,3){\line(-1,0){1.5}}
\put(11,2){\line(0,1){0.5}}
\put(3,2.5){\line(0,1){0.5}}
\put(1.5,3){\line(0,1){0.5}}
\put(1.5,3.5){\line(-1,0){1.3}}

\put(0.1,-0.2){\textbf{Figure 2.1}: The $E_{\infty}$-page of Greenlees spectral sequence for $ko_{*}(BSD_{16}).$}
\end{picture}
\end{center}

For the precise generators description of $[n]$ and $2^{r}$ in Figure 2.1, we refer the reader to Lemma 5.4.3 in \cite{kijti}.  Now, the results can be read from the $E_{\infty}$-page directly with the filtration given by
$$ko_{n}(BSD_{16})=F^{n}_{0}\supseteq F^{n}_{1}\supseteq F^{n}_{2}\supseteq F^{n}_{3}=0, $$ with
$ F^{n}_{0}/F^{n}_{1}\cong E^{0,n}_{\infty}, F^{n}_{1}/F^{n}_{2}\cong E^{-1,n+1}_{\infty}$ and $F^{n}_{2}\cong E^{-2,n+2}_{\infty}$.  Precisely, we use two short exact sequences to determine $ko_{n}(BSD_{16})$, viz;
$$ 0\longrightarrow F^{n}_{1} \longrightarrow ko_{n}(BSD_{16})\longrightarrow E^{0,n}_{\infty}\longrightarrow0, $$
and
$$ 0\longrightarrow E^{-2,n+2}_{\infty} \longrightarrow F^{n}_{1}\longrightarrow E^{-1,n+1}_{\infty}\longrightarrow0. $$
From this fact and $E_{\infty}$-page, we see that there are extension problems in degree $n\geq8$ when $n$ is congruent to $0,1,2$ modulo $8$.  It can be shown that they are all trivial, (Section 5.5.2, 6.5 and Lemma 6.5.1 in \cite{kijti}).  This implies that $ko_{*}(BSD_{16})\cong E_{\infty}$-page and hence by lemma \ref{lem Ap},
\begin{equation}\label{kerAp}
    \ker(Ap)\cong H^{1}_{I}(ko^{*}(BSD_{16}))\oplus H^{2}_{I}(ko^{*}(BSD_{16})).
\end{equation}

Explicitly, we have (cf. \cite{kijti}, Theorem 6.5.2)
\begin{center}
\begin{tabular}{|c|c|c|}
  \hline
   & & \\
  $n$ & $[H_{I}^{1}(ko^{*}(BSD_{16}))]_{n+1}$ & $[H_{I}^{2}(ko^{*}(BSD_{16}))]_{n+2}$ \\
  & & \\
  \hline
  \hline
   $n\leq 3$ & $0$ & $0$ \\
  $3$ & $[4]\oplus[8]\oplus[8]$ & $0$ \\
  $4$ & $0$ & $2$ \\
  $5$ & $[2]$ & $0$ \\
  $6$ & $0$ & $0$ \\
  $7$ & $[2]\oplus[4]\oplus[16]\oplus[32]$ & $0$ \\
  $8$ & $2$ & $2$ \\
  $9$ & $2\oplus [2]$ & $0$ \\
  $10$& $0$ & $2$ \\
  $11$& $[8]\oplus[16]\oplus[128]\oplus[128]$ & $0$ \\
  $12$& $0$ & $2^{2}$ \\
  $13$ & $[4]$ & $0$ \\
  $14$ & $0$ & $2$ \\
  $15$ & $[2]\oplus[8]\oplus[16]\oplus[256]\oplus[512]$ & $0$ \\
  $8k\geq16$ & $2^{2}=2 \oplus \eta(ko_{8k-1}(BSD_{16}))$ & $2^{k}$ \\
  $8k+1\geq 17$ & $2^{2} \oplus [2^{k}]=\eta(\widetilde{ko}_{8k}(BSD_{16}))\oplus [2^{k}]$ & $0$ \\
  $8k+2\geq 18$ & $0$ & $2^{k}$ \\
  $8k+3\geq 19$ & $[2^{k-1}]\oplus[2\cdot 4^{k}]\oplus[4^{k+1}]\oplus[8\cdot16^{k}]\oplus[8\cdot16^{k}]$ & $0$ \\
  $8k+4\geq 20$ & $0$ & $2^{k+1}$ \\
  $8k+5\geq 21$ & $[2^{k+1}]$ & $0$ \\
  $8k+6\geq 22$ & $0$ & $2^{k}$ \\
  $8k+7\geq 23$ & $[2^{k}]\oplus[2\cdot 4^{k}]\oplus[4^{k+1}]\oplus[16^{k+1}]\oplus[2\cdot16^{k+1}]$ & $0$ \\
  \hline
  \hline
\end{tabular}.
\end{center}
\begin{center}
    \textbf{Table 2.2:} The additive structure of $\ker(Ap)$ for $G=SD_{16}$.
\end{center}

Here, $[H_{I}^{1}(ko^{*}(BSD_{16}))]_{n+1}$ and $[H_{I}^{2}(ko^{*}(BSD_{16}))]_{n+2}$ contribute to $(\ker{Ap})_{n}\subseteq ko_{n}(BSD_{16})$.  Note further from \cite{kijti} (Theorem 6.5.2) that the 2-column, $[H_{I}^{2}(ko^{*}(BSD_{16}))]_{n+2}$, is embedded in $H_{n}(BSD_{16};\mathbb{F}_{2})$.  For the one column, $[H_{I}^{1}(ko^{*}(BSD_{16}))]_{n+1}$, we have that the generators of $[H_{I}^{1}(ko^{*}(BSD_{16}))]_{8+1}$ is embedded in $H_{8}(BSD_{16};\mathbb{F}_{2})$, namely $\xi(yuP):=(yuP)^{\vee}$ (the dual element of $yuP \in H^{8}(BSD_{16};\mathbb{F}_{2})$)\footnote{See $H^{*}(BSD_{N};\mathbb{F}_{2})$ in Proposition \ref{HSDN}.}.  In fact, this generator is an element in $H_{I}^{1}(\tau)$, (see (\ref{local cohomology of sd16}) above), by the GSS calculation; precisely, it is $\widetilde{\eta}[\overline{u}_{4}]/(\overline{u}_{4})^{2}$, Lemma 5.2.4 in \cite{kijti}.  But the BSS calculation via $ku_{*}(BSD_{16})$, Lemma 6.3.1 and Figure 6.3 in \cite{kijti}, reveals that it is embedded in $H_{n}(BSD_{16};\mathbb{F}_{2})$.

In higher degrees, $[H_{I}^{1}(ko^{*}(BSD_{16}))]_{(8k)+1}$ contains two generators of order 2 where one of them is an $\eta$-multiple (Theorem 6.5.2 in \cite{kijti}) and the other one is embedded in $H_{8k}(BSD_{16};\mathbb{F}_{2})$ as $\xi(yuP^{2k-1}):=(yuP^{2k-1})^{\vee}$, (analogous reference as the case $\xi(yuP):=(yuP)^{\vee}$).  The generator of $[H_{I}^{1}(ko^{*}(BSD_{16}))]_{9+1}$ is an $\eta$-multiple and in higher degree $[H_{I}^{1}(ko^{*}(BSD_{16}))]_{(8k+1)+1}$ contains one $\eta$-multiple generator and one $\eta^{2}$-multiple generator (coming from the $\eta$-multiple in degree $(8k)+1$).  The other generators in this column which we did not mention above are detected by $H_{I}^{1}(\overline{QO})$ which can be dealt with by character theory, see more details about the generators in Theorem 6.5.2 in \cite{kijti}.

Note also that the part $TO$ detected in ordinary cohomology of $SD_{2^{N}}$ stays the same for all $N\geq4$ because the generators and relations of $H^{*}(BSD_{2^{N}};\mathbb{F}_{2})$ do not depend on $N$.  Indeed, the $2-$column of the $E_{\infty}$ page of the Greenlees spectral sequence will be unchanged for all semi-dihedral groups, since we can realize all of them, as in the case of $SD_{16}$, by positive scalar metric manifolds (see Section 5) and thus no more non-zero differentials of the GSS detects them for general $N$. However, the representation theory for the calculation of $H_{I}^{1}(\overline{QO})$ for semidihedral groups of order $2^{N}$, $N> 4$, is increasingly complicated.

\section{The one column and lens spaces detected in periodic K-theory}
In order to prove the conjecture for groups with periodic cohomology, Botvinnik, Gilkey and Stolz \cite{2} used the eta invariant to detect the orders of fundamental classes of lens spaces and lens space bundles. Since the 1-column contains the part of $Ker(Ap)$ detected in periodic K-theory, by restricting representations to periodic subgroups, we can use the naturality of the eta invariant to calculate the orders of the lens spaces which we use for the periodic subgroups, under induction.

We have that $SD_{16}=<s,t;s^8=t^2=1, tst=s^3>$. First consider the cyclic subgroup $<t>$ of order $2$. We know that the real projective space $\RP^{4k+3}$ has fundamental group cyclic of order two, and thus there is a classifying map $u:\RP^{4k+3} \rightarrow B<t>$. We can compose this with the inclusion into $SD_{16}$ to get a map $\RP^{4k+3} \rightarrow B<t> \hookrightarrow BSD_{16}$, and then use the eta invariant to ask how big a subgroup we span by the inclusion of this fundamental class.

Similarly, we can take the cyclic subgroup $<s>$ of order $8$. Lens spaces of the form $S^{4k+3}/C_8$ together with lens space bundles of the form $S^{4k+3}/C_8 \rightarrow M^{4k+5} \rightarrow S^2$ map naturally into $B<s>$, and we can again compose with inclusion and use the eta invariant formulae of \cite{2}. Finally we can also the quaternion subgroup $<s^2,ts>$ of order $8$, for which we have lens spaces of the form $S^{4k+3}/Q_8$ mapping into the classifying space. We will show that this collection of manifolds realizes all of $Ker(A) \cap Im(p)$.

We start by recalling the basic properties and relevant formulae for the eta invariant.
Atiyah, Patoudi and Singer showed that there is a formula for the index of the Dirac operator $D$ for a manifold $W$ with boundary $M$, analogous to the usual index formula, but with a correction term $\eta(D(M))$, known as the eta invariant, depending only on the boundary.

Given a representation $\rho$ of a discrete group $\pi$ and a map from a manifold $f:M \to B\pi$, we can form a vector bundle $V_{\rho}:\widetilde{M}\times_{\pi} \rho \to M$, where $\widetilde{M}$ is the $\pi$ cover of $M$ classified by $f$. We can then consider the Dirac operator $D(M,f,\rho)$, which is the Dirac operator of $M$ twisted by $V_{\rho}$, and its Eta invariant $\eta(M,f)(\rho)$.

We recall from \cite{7} that there is a geometric description of periodic K-theory as follows, where the isomorphism is induced by $pD$:
$$\Omega_{*}^{spin}(X)/T_{*}(X)[B^{-1}] \cong KO_{*}(X)$$
where $T_{*}(X)$ is the subgroup generated by quaternion projective bundles, and $B=B^8$ is the Bott manifold, which is any simply connected manifold with $\hat{A}(B)=1$.
Using this description, we have the following result \cite{2}:
\begin{thm}\label{thm of eta Z and Zmod2}
Let $\rho$ be a virtual representation of $\pi$ of virtual dimension zero. Then the homomorphisms
$$\eta(\rho):\Omega_{n}^{spin}(B\pi)\rightarrow \R/\Z;\eta(\rho):KO_{n}(B\pi)\rightarrow \R/\Z$$
sending $[f:M \rightarrow B\pi]$ to $\eta(M,f)(\rho)$ are well defined. Further if $\rho$ is real and $n \equiv 3 \mod 8$, or if $\rho$ is quaternion and $n \equiv 7 \mod 8$, then we can replace the range of $\eta(\rho)$ by $\R/2\Z$, meaning that the map $\eta(\rho):\Omega_{n}^{spin}(B\pi)\rightarrow \R/2\Z$ is well defined, and similarly for $KO$.
\end{thm}
The following theorem, deduced from work of Donnelly \cite{5}, is the main tool for actually computing some eta invariants.
\begin{thm}\label{eta calculation thm}
Let $\rho$, $\pi$ be as above, and $\tau:\pi \rightarrow U(m)$, be a fixed point free representation. Assume there exists a representation $det(\tau)^{1/2}$ of $\pi$ whose tensor square is $det(\tau)$. Then letting $M=S^{2n-1}/\tau(\pi)$ with the inherited structures, we have
$$\eta(M)(\rho)=\mid \pi \mid^{-1} \sum_{1 \neq g \in \pi} \frac{Trace(\rho(g))det(\tau(g))^{1/2}}{(det(I-\tau(g)))}$$
\end{thm}

Botvinnik, Gilkey and Stolz \cite{2} use this result to prove the conjecture for generalized quaternion groups, and we shall now outline a modified version of their proof for $Q_8$, so as to understand explicitly the order of the subgroup of $ko_{4k+3}(BSD_{16})$ obtained by restricting representations. Here is the character table of $Q_8$

\begin{center}
\begin{tabular}{|c|c|c|c|c|c|}
\hline
$ $&1 & 1 & $2$ & $2$ & $2$ \\
$\rho$&$1$ & $-1$ & $[i]$ & $[j]$ & $[k]$\\
\hline
$1=\rho_0$ & $1$ & $1$ & $1$ & $1$ & $1$\\
$\kappa_1$ & $1$ & $1$ & $-1$ & $1$ & $-1$\\
$\kappa_2$ & $1$ & $1$ & $1$ & $-1$ & $-1$\\
$\kappa_3$ & $1$ & $1$ & $-1$ & $-1$ & $1$\\
$\tau$ & $2$ & $-2$&$0$ &$0$& $0$\\
\hline
\end{tabular}
\end{center}

We have three cyclic subgroups $H_t$, with $t=1,2,3$ generated by $i,j,k$ respectively, and thus by viewing $S^{4m-1}$ inside $\qt^m$ we can consider the manifolds $M_{t}^{4m-1}=S^{4m-1}/H_t$ together with $M_{Q}^{4m-1}=S^{4m-1}/Q_8$. We define quadruples of eta invariants as follows:
$$ \overrightarrow{\eta}(M)=(\eta(M)(\rho_0-\kappa_1),\eta(M)(\rho_0-\kappa_3),\eta(M)(2-\tau),\eta(M)(2-\tau)^2)$$
Since the representations $(2-\tau)^{2a}$ are real and $(2-\tau)^{2a+1}$ are quaternion, theorem \ref{thm of eta Z and Zmod2} tells us that if $n \equiv 3 \mod 8$
$$\overrightarrow{\eta}(M^n) \in \R/2\Z \oplus \R/2\Z \oplus \R/\Z \oplus \R/2\Z $$
while in $n \equiv 7 \mod 8$
$$\overrightarrow{\eta}(M^n) \in \R/\Z \oplus \R/\Z \oplus \R/2\Z \oplus \R/\Z $$

We now proceed a little differently now than in \cite{2}, because it is useful to resolve the extensions for $Q_8$, so as to understand the inclusion better. This is done in the thesis of Bayen \cite{1} using the Adams spectral sequence, and we now outline the results, and verify his calculations from a geometric viewpoint.

There is a 2-local stable decomposition of $BQ_8$, and higher groups, into indecomposable summands as follows \cite{3}:
$$BQ_l=BSL_2(q) \vee 2\Sigma^{-1}BS^3/BN$$
where $N$ is the normalizer of a maximal torus in $S^3$, and $l$ is the largest power of $2$ dividing $q^2-1$, for $q$ an odd prime power. Note that $q=3$ for $l=8$. The method then is to make the calculation for each summand. It then turns out that for $k=3,7$ we get that 2-locally:
$$ko_{8m+k}(\Sigma^{-1}BS^3/BN)=\Z_{2^{2m+2}}$$
$$ko_{8m+3}(BSL_2(3))=\Z_{2^{4m+3}}\oplus \Z_{2^{2m}}=<x_{8m+3}> \oplus <\beta x_{8m-5}-16x_{8m+3}>$$
$$ko_{8m+7}(BSL_2(3))=\Z_{2^{4m+6}}\oplus \Z_{2^{2m}}=<x_{8m+7}> \oplus <\beta x_{8m-1}-16x_{8m+7}>$$
Then it follows $ko_*(BQ_8)=ko_*(BSL_2(3))\oplus 2 ko_*(\Sigma^{-1}BS^3/BN)$.

We shall now give a proof of the GLR conjecture for $Q_8$, by describing the generators as images of fundamental classes of manifolds.

Firstly, in the splitting $BQ_l=BSL_2(q) \vee 2\Sigma^{-1}BS^3/BN$, the $2\Sigma^{-1}BS^3/BN$ term is independent of $l$, and thus so is $2 ko_*(\Sigma^{-1}BS^3/BN)$ (here $ko_*$ is 2-locally), and for $n=8m+3,8m+7$ we have the following calculation from \cite{2}:
$$\overrightarrow{\eta}(M_{1}^n-M_{2}^n)=(2^{-2m-l},0,0,0)$$
$$\overrightarrow{\eta}(M_{1}^n-M_{3}^n)=(*,2^{-2m-l},0,0)$$
Where * is a term we aren't interested in. Here $l=1,2$ for $n=8m+3,8m+7$ respectively. Since $\rho_0-\kappa_1$ and $\rho_0-\kappa_3$ are both real however, the subspace spanned has order $(2^{2m+2})^2$ in both cases. Thus it follows that $2$-locally, the manifolds $M_{1}^n-M_{2}^n, M_{1}^n-M_{3}^n$ generate all of $2 ko_n(\Sigma^{-1}BS^3/BN)$, with $ n \equiv 3\mod 4$.

It thus remains to realize the $BSL_2(3)$ part by some fundamental classes of manifolds, and so we consider the manifolds $M_Q^{4k+3}=S^{4k+3}/Q_l$. We claim that the remaining part of $ko_{4k+3}(BQ_l) \cong ko_{4k+3}(BSL_2(3))$ can be realized by the images of the fundamental classes of the manifolds $M_Q^{4k+3}$ and $B^8 \times M_Q^{4k-5}$. We do this by explicit eta calculations, using only representations $(2-\tau)^a$. Note that this suffices since from above and \cite{2}, $\eta(M_{1}^n-M_{i}^n)(2-\tau)^a=0$ for any $a$. 

We start by calculating the order of $[M_Q^{4k+3}]$. Note that $M_Q^{4k+3}=S^{4k+3}/(k+1)\tau$, where $(k+1)\tau$ is the $k+1$-fold sum $\tau \oplus \cdots \oplus \tau$. Now $\tau$ is a unitary representation which may be characterized by 
\begin{displaymath}
\tau(i)=\left(\begin{array}{ccc}
i & 0 \\
0 & -i \\
\end{array} \right), \quad
\tau(j)=\left(\begin{array}{ccc}
0 & \omega^3 \\
\omega & 0 \\
\end{array} \right),
\end{displaymath}
where $\omega=(1+i)/\sqrt{2}$. Thus for $x \neq 1$ we get $det(1-\tau)(x)=4$ if $x=-1$, and $2$ else. So we can apply Theorem \ref{eta calculation thm} to see:
$$\begin{array}{rcl}
  \eta(M_Q^{4k+3})(2- \tau) & = & \displaystyle\frac{1}{8}\sum_{1 \neq g \in Q_8} \frac{Trace((2-\tau)(g))det((k+1)\tau(g))^{1/2}} {(det(I-(k+1)\tau(g)))} \\
   & = & (1/8)(4/4^{k+1}+6(2/2^{k+1}))=1/2^{2k+3}+3/2^{k+2}
\end{array}$$
which has order $2^{2k+3} \in \R/\Z$, and $2^{2k+4} \in \R/2\Z$. Thus in dimensions $8m+3$ the $ko-$fundamental class $[M_Q^{8m+3}]$ has order at least $2^{4m+3} \in ko_*(BQ_8)$, and the calculations in \cite{3} imply that it has exactly this order. The representation $\tau$ is quaternion, so in dimensions $8m+7$ the $ko-$fundamental class $[M_Q^{8m+7}]$ has order $2^{4m+6}$, since the eta invariant extends to $\R/2\Z$.

We wish to see how large a subspace is spanned by $M_Q^{4k+3}$ and $B \times M_Q^{4k-5}$. Multiplying by the Bott element is a monomorphism, so $B \times M_Q^{4k-5}$ has the same order as $M_Q^{4k-5}$. Then, just as above, we can make the following calculations, remembering $(2-\tau)^a$ is real for $a$ even, and quaternion for $a$ odd:
$$\begin{array}{rcl}
\eta(M_Q^{4k+3})(2- \tau)^2&=&\displaystyle\frac{1}{8}\sum_{1 \neq g \in Q_8} \frac{Trace((2-\tau)(g))^2det((k+1)\tau(g))^{1/2}} {(det(I-(k+1)\tau(g)))} \\
&=&(1/8)(16/4^{k+1}+6(4/2^{k+1}))=2/4^{k+1}+3/2^{k+1} \in \R/2\Z, \end{array}$$

$$\begin{array}{rcl}
\eta(M_Q^{4k+3})(2- \tau)^3&=&\displaystyle\frac{1}{8}\sum_{1 \neq g \in Q_8} \frac{Trace((2-\tau)(g))^3det((k+1)\tau(g))^{1/2}} {(det(I-(k+1)\tau(g)))} \\
&=&(1/8)(64/4^{k+1}+6(8/2^{k+1}))=2/4^{k}+6/2^{k+1} \in \R/2\Z. \end{array}$$
Thus using these and the earlier formula we can say:
$$(\eta(M_Q^{8m+3})(2- \tau),\eta(M_Q^{8m-5} \times B^8)(2- \tau))=(1/2^{4m+3}+3/2^{2m+2},1/2^{4m-1}+3/2^{2m}),$$
while for $(2-\tau)^2$ we get:
$$(\eta(M_Q^{8m+3})(2- \tau)^2,\eta(M_Q^{8m-5} \times B^8)(2- \tau)^2)=(1/2^{4m+2}+3/2^{2m+2},1/2^{4m-2}+3/2^{2m}).$$
The order of the subgroup spanned may be bounded below by the order $\in \R/\Z$ of the determinant of the following matrix in $\R/\Z$, since the two rows are the images of the two homomorphisms $\eta(2-\tau), \eta(2-\tau)^2$:
\begin{displaymath}
X=\left(\begin{array}{ccc}
\eta(M_Q^{4k+3})(2- \tau) & \eta(M_Q^{4k-5} \times B^8)(2- \tau) \\
\eta(M_Q^{4k+3})(2- \tau)^2 & \eta(M_Q^{4k-5} \times B^8)(2- \tau)^2 \\
\end{array} \right)
\end{displaymath}
From our calculations we can read this matrix off. When $4k+3=8m+3$ we get:\\
\begin{displaymath}
X=\left(\begin{array}{ccc}
1/2^{4m+3}+3/2^{2m+2} & 1/2^{4m-1}+3/2^{2m} \\
1/2^{4m+2}+3/2^{2m+2} & 1/2^{4m-2}+3/2^{2m}\\
\end{array} \right)
\end{displaymath}
Now the determinant of $X$ is $$1/2^{8m+1}+3/2^{6m+3}+3/2^{6m}+9/2^{4m+2}-1/2^{8m+1}-3/2^{6m+2}-3/2^{6m+1}-9/2^{4m+2}$$
$$=3/2^{6m+3} + \cdots$$
The dots mean terms of lower order, so the determinant has order $2^{6m+3} \in \R/\Z$ just as required.

Analogously in dimensions $8m+7$ we have:
$$(\eta(M_Q^{8m+7})(2- \tau),\eta(M_Q^{8m-1} \times B^8)(2- \tau))=(1/2^{4m+6}+3/2^{2m+4},1/2^{4m+2}+3/2^{2m+2})$$
while for $(2-\tau)^2$ we get:
$$(\eta(M_Q^{8m+7})(2- \tau)^3,\eta(M_Q^{8m-1} \times B^8)(2- \tau)^3)=(1/2^{4m+3}+3/2^{2m+2},1/2^{4m-1}+3/2^{2m})$$
The analogous determinant then has order $2^{6m+6}$ in $\R/\Z$ as required. Thus $2-$locally, all of $ko_{4k+3}(BQ_8)/2 ko_{4k+3}(\Sigma^{-1}BS^3/BN)$ is realized by the manifolds $M_Q^{4k+3},M_Q^{4k-5} \times B^8$, and may be detected by the eta invariants, using only the virtual representations $2-\tau$ and $(2-\tau)^2$.

Now consider cyclic groups of order $l=2^k$, which we identify with the subgroup of $S^1$ consisting of $l$-th roots of unity:
$$C_l=\{\lambda \in S^1 | \lambda^l=1\}$$
For an integer $a$ we let $\rho_a$ be the representation of $S^1$ where $\lambda \in S^1$ acts by multiplication by $\lambda^a$. For a tuple of integers $\overrightarrow{a}=(a_1, \cdots ,a_t)$, the representation $\lambda^{a_1 }\oplus \cdots \oplus \lambda^{a_t}$ restricts to a free $C_l$ action on $S^{2t-1}$ if and only if all the $a_j$ are odd. Let $t=2i$ be even, and define the quotient manifold
$$X^{4i-1}(l,\overrightarrow{a})=S^{2t-1}/(\rho_{a_1} \oplus \cdots \oplus \rho_{a_t})(C_l)$$
This is a lens space with a natural positive scalar curvature metric, and it inherits a natural spin structure in dimension $4i-1$.

Note that in dimensions $4k+1$ the corresponding construction does not yield a spin manifold. So, consider the vector bundle $H \otimes H \oplus (2k-1)\C \rightarrow S^2$, where $H$ is the Hopf line bundle. As above, for a tuple of integers $\overrightarrow{a}=(a_1, \cdots ,a_t)$, the representation $\lambda^{a_1 }\oplus \cdots \oplus \lambda^{a_t}$ restricts to a free $C_l$ action on the sphere bundle $S(H \otimes H \oplus (2k-1) \C \rightarrow S^2)$ if and only if all the $a_j$ are odd. The quotient
$$X^{4i+1}(l,\overrightarrow{a})=S(H \otimes H \oplus (2k-1)\C \rightarrow S^2)/C_l$$
is a positive scalar curvature spin manifold of dimension $4k+1$ which is a lens space bundle over the two-sphere. A general formula for the eta invariants of manifolds of the form $S(H_1 \oplus \cdots \oplus H_k \rightarrow S^2)/C_l$, where the $H_i$ are complex line bundles and the $C_l$ action is specified by a tuple of integers $\overrightarrow{a}$ as above, is given by \cite{2}:
\begin{thm}
Let $\rho \in R_0(C_l)$ and $M=S(H_1 \oplus \cdots H_k \rightarrow S^2)/C_l(\overrightarrow{a})$ as above. Then
$$\eta(M)(\rho)=l^{-1}\sum_{1 \neq \lambda \in C_l} Trace(\rho(\lambda) \frac{\lambda^{(a_1+ \cdots a_{2i})/2}}{(1-\lambda^{a_1})\cdots (1-\lambda^{a_{2i}})}. \sum_j \frac{1}{2}c_1(H_j)[\CP_1] \frac{1+\lambda^{a_j}}{1-\lambda^{a_j}}$$
\end{thm}

 It is then clear that $\pi_1(X^{4i-1}(l,\overrightarrow{a}))=C_l$, giving these manifolds natural $C_l$ structures, and in order to apply theorem $3.2$ we can simply define $det(\rho_{a_1} \oplus \cdots \oplus \rho_{a_t})^{1/2} :=\rho_{(a_1 + \cdots a_{2i})/2}$. We say $L^{n}=X^{n}-X_{0}^{n} \in \Omega_{n}^{spin}(BC_l)$, where $X_{0}^{n}$ is the same manifold with trivial $C_l$ structure. We can then use the following formula: 

 For $\lambda \in C_l$ we define $$f_{4i-1}(\overrightarrow{a_{2i}})(\lambda):=\lambda^{(a_1 + \cdots a_{2i})/2}(1-\lambda^{a_1})^{-1}\cdots(1-\lambda^{a_{2i}})^{-1}$$ and $f_{4i+1}(\overrightarrow{a_{2i}})(\lambda):=f_{4i-1}(\overrightarrow{a_{2i}})(\lambda)(1+\lambda^{a_1})(1-\lambda^{a_1})^{-1}$. Then we have from the above two theorems:
\begin{lem}
$\eta(L^{n}(l,\overrightarrow{b_{2i}}))(\rho)=\displaystyle l^{-1}\sum_{1 \neq \lambda \in C_l}f_{4i \pm 1}(\overrightarrow{b_{2i}})(\lambda)Trace(\rho(\lambda))$
\end{lem}

We remark that the eta invariant is natural with respect to inclusions \cite{1}, \cite{2}, meaning that if we have an inclusion of groups $ f:H\hookrightarrow G$, and a class represented by a manifold $[M] \in \Omega_{n}^{spin}(BH)$, then for a representation $\rho$ of $G$ the pull-back bundle over $M$ is determined by restricting the representation and we have $\eta([M])(f^*(\rho))=\eta(f_*([M]))(\rho)$, as should be clear from the following pull-back diagram:
$$
\xymatrix{
i^*(EH \times_H f^*(\rho)) \ar[r] \ar [d] & EH \times_H f^*(\rho) \ar[r] \ar [d] & EG \times_G \rho \ar[d] \\
M \ar[r]^{i} & BH \ar[r]^{f} & BG \\
}
$$
We are now ready to prove the main result of this section, which verifies the conjecture in odd dimensions. We start with the character table of $SD_{16}=<s,t;s^8=t^2=1,tst=s^3>$.
\begin{center}
\begin{tabular}{|c|c|c|c|c|c|c|c|}
\hline
$ $& $1$ & $1$ & $2$ & $2$ & $2$ & $4$ & $4$\\
$\rho$&$[1]$ & $[s^4]$ & $[s]$ & $[s^2]$ & $[s^5]$ & $[t]$ & $[ts]$\\
\hline
\hline
$1=\rho_0$ & $1$ & $1$ & $1$ & $1$ & $1$& $1$ & $1$\\
$\chi_2=\hat{C_8}$ & $1$ & $1$ & $1$ & $1$ & $1$ & $-1$ & $-1$\\
$\chi_3=\hat{D_8}$ & $1$ & $1$ & $-1$ & $1$ & $-1$ & $1$ & $-1$\\
$\chi_4=\hat{Q_8}$ & $1$ & $1$ & $-1$ & $1$ & $-1$ & $-1$ & $1$\\
$\chi_{\rho}$ & $2$ & $-2$ & $\sqrt{2}i$ & $0$ & $-\sqrt{2}i$ & $0$ & $0$\\
$\chi_{\rho^2}$ & $2$ & $2$ & $0$ & $-2$ & $0$ & $0$ & $0$\\
$\chi_{\rho^5}$ & $2$ & $-2$ & $-\sqrt{2}i$ & $0$ & $\sqrt{2}i$ & $0$ & $0$\\
\hline
\hline
\end{tabular}
\end{center}

\begin{prop}
Inclusion of $ko-$fundamental classes of the lens spaces given above from the cyclic subgroups $<s>,<t>$ and the quaternion subgroup $<s^2,ts>$, spans all of $ko_{4m+3}(BSD_{16})$. Inclusion from $<s>$ spans all of $ko_{8m+5}(BSD_{16})$, and together with eta multiples, all of $Ker(Ap) \subset ko_{8m+9}(BSD_{16})$
\end{prop}
\begin{proof}
We recall that $|ko_{8m+3}(BSD_{16})|=2^{8+13m}$, while $|ko_{8m+7}(BSD_{16})|=2^{12+13m}$.\\
Start by considering explicitly the restrictions to $Q_8$, and using the calculations known for $Q_8$. Recall that we had
$$\overrightarrow{\eta}(M)=(\eta(M)(\rho_0-\kappa_1),\eta(M)(\rho_0-\kappa_3),\eta(M)(2-\tau),\eta(M)(2-\tau)^2)$$
The real representation $\chi_{\rho^2}$ restricts to $\kappa_1+\kappa_3$, so since $\overrightarrow{\eta}(M_{1}^n-M_{2}^n)=(2^{-2m-1}/2^{-2m-2},0,0,0)$ in $n=8m+3$, this evaluates to $2^{-2m-1} \in \R/2\Z$ which has order $2^{2m+2}$, while in $n=8m+7$ we get $2^{-2m-2} \in \R/\Z$ which again has order $2^{2m+2}$ .

Further, the representations $\chi_{\rho}, \chi_{\rho^5}$ both restrict to the natural representation $\tau$ of $Q_8$, and $(2-\tau)^2=4+\tau^2-4\tau$ is restricted to by $4+\chi_{\rho}.\chi_{\rho^5}-2(\chi_{\rho} + \chi_{\rho^5})$ which is a real representation since $\chi_{\rho}, \chi_{\rho^5}$ are complex conjugate.
So, by restricting representations we see directly that in $n=8m+3$ we have spanned a subspace of order $2^{2m+2}8^{2m+1}$ by including from $Q_8$.

In $n=8m+7$, since $\chi_{\rho}, \chi_{\rho^5}$ are not quaternion representations, we must halve the order of whatever is in $\R/2\Z$ from the $2-\tau$, giving a subspace of order $2^{2m+2}8^{2m+2}2^{-1}=2^{2m+1}8^{2m+2}$.

Note further that the representation $\chi _4$ of $SD_{16}$ restricts trivially on the quaternion subgroup $<s^2,ts>$ and to the non trivial representation on $<t>$, and thus in addition to what is induced from $Q_8$ we have $[\RP^n] \rightarrow ko_n(B<t>) \hookrightarrow ko_n(BSD_{16})$ of orders $2^{4m+3}, 2^{4m+4}$ when $n=8m+3, 8m+7$ respectively.

Putting this all together, we have in $8m+3$ a subgroup of order $2^{2m+2}8^{2m+1}2^{4m+3}=2^{8+12m}$ (this is already enough in $n=3$), while in $8m+7$ the order is $2^{4m+4} 2^{2m+1}8^{2m+2}=2^{11+12m}$. We now realize what is remaining by using  the lens spaces $L^{n}$ described above, viewing the fundamental group $C_8=<s>$ sitting inside $SD_{16}$.

We start by calculating the eta invariant with respect to the real representation $\rho_0-\rho_4$ of $C_8$. We show that $\eta(L^{8m+j}(l=8,\overrightarrow{a}))(\rho_{4}-\rho_0)$ has order at least $2^m$ for both $j=3,7$ for suitable $\overrightarrow{a}$. Firstly, let $\omega=(1+i)/\sqrt{2}$ be the generator of $C_8$. Then $\rho_4(\omega)=-1, \omega^5=-\omega, \omega^7=-\omega^3$ and we have:
$$\eta(L^{3}(l=8,(1,1)))(\rho_{4}-\rho_0)=\frac{1}{8}\{\frac{2\omega}{(1-\omega)^2}+\frac{2\omega^3}{(1-\omega^3)^2}+\frac{2\omega^5}{(1-\omega^5)^2}+\frac{2\omega^7}{(1-\omega^7)^2}\}$$
$$=\frac{1}{4}\{\frac{\omega((1+\omega)^2-(1-\omega)^2)}{(1-i)^2}+\frac{\omega^3((1+\omega^3)^2-(1-\omega^3)^2)}{(1+i)^2}\}$$
$$=\frac{1}{8i}\{-\omega(4\omega)+\omega^3(4\omega^3)\}=-8i/8i=-1$$
which has order $2 \in \R/2\Z$. Similarly, in $n=7$ we have:
$$\eta(L^{7}(l=8,(1,1,1,1)))(\rho_{4}-\rho_0)=\frac{1}{8}\{\frac{2\omega^2}{(1-\omega)^4}+\frac{2\omega^6}{(1-\omega^3)^4}+\frac{2\omega^{10}}{(1-\omega^5)^4}+\frac{2\omega^{14}}{(1-\omega^7)^4}\}$$
$$=\frac{1}{4}\{\frac{i((1+\omega)^4+(1-\omega)^4)}{(1-i)^4}-\frac{i((1+\omega^3)^4+(1-\omega^3)^4)}{(1+i)^4}\}$$
$$=\frac{-1}{16}\{12i^2-12i(-i)\}=24/16=3/2$$
which has order $2 \in \R/\Z$.

For the general case we let $\overrightarrow{a}=(a_1, \cdots,a_{2k}), K=(a_1+\cdots a_{2i})/2$ and use the following inductive trick:
$$\eta(L^{4k+7}(l=8,(\overrightarrow{a},1,1,5,5)))(\rho_{4}-\rho_0)=$$
$$=\frac{1}{8}\{\frac{2\omega^{K+6}}{(1-\omega^{a_1})\cdots(1-\omega^{a_{2k}})(1-i)^2}+\frac{2\omega^{3K+18}}{(1-\omega^{3a_1})\cdots(1-\omega^{3a_{2k}})(1+i)^2}+$$
$$\frac{2\omega^{5K+30}}{(1-\omega^{5a_1})\cdots(1-\omega^{5a_{2k}})(1-i)^2}+\frac{2\omega^{7K+42}}{(1-\omega^{7a_1})\cdots(1-\omega^{7a_{2k}})(1+i)^2}\}$$
Noting again that $(1 \pm i)^2=\pm 2i$, and $\omega^{4j+2}=i,-i$ for $j$ even, odd respectively, we may take the common terms out of the bracket to simplify:
$$=\frac{i}{16i}\{\frac{2\omega^{K}}{(1-\omega^{a_1})\cdots(1-\omega^{a_{2k}})}+\frac{2\omega^{3K}}{(1-\omega^{3a_1})\cdots(1-\omega^{3a_{2k}})}+$$
$$\frac{2\omega^{5K}}{(1-\omega^{5a_1})\cdots(1-\omega^{5a_{2k}})}+\frac{2\omega^{7K}}{(1-\omega^{7a_1})\cdots(1-\omega^{7a_{2k}})}\}$$
$$=(1/2) \eta(L^{4k-1}(l=8,(\overrightarrow{a}))(\rho_{4}-\rho_0)$$
Thus $\eta(L^{4k+7}(l=8,(\overrightarrow{a},1,1,5,5)))(\rho_{4}-\rho_0)$ has twice the order in $\R/\Z$ of $\eta(L^{4k-1}(l=8,(\overrightarrow{a}))(\rho_{4}-\rho_0)$, and the claim is immediate by induction.

Thus we can consider $6-$tuples of eta invariants as follows. The three non-trivial one dimensional representations of $SD_{16}$ can each be characterized by their kernels, which are $C_8, D_8$ and $Q_8$ respectively. We will denote the representation with kernel $C_8$ as $\hat{C_8}$ for example. Then for a manifold $M$ we set:
$$ \overrightarrow{\eta}(M)=(\eta(M)(1-\hat{D_8}),\eta(M)(1-\hat{C_8}),\eta(M)(1-\hat{Q_8}),
\eta(M)(2-\chi_{\rho}),\eta(M)(4+\chi_{\rho}.\chi_{\rho^5}-2(\chi_{\rho} + \chi_{\rho^5}))$$
Then using \cite{2} and the calculations we've already given we have upto order at least in $n=8m+3$:
$$\begin{array}{ccl}
\overrightarrow{\eta}(L^{8m+3},\overrightarrow{a})&=&(2^{-m-1},0,2^{-m-1},*,*,*) \\
\overrightarrow{\eta}(\RP^n)&=&(0,2^{-4m-3},2^{-4m-3},2^{-4m-3},2^{-4m-2},2^{-4m-3})\\
\overrightarrow{\eta}(M_{1}^n-M_{2}^n)&=&(*,*,0,2^{-2m-2},0,0)\\
\overrightarrow{\eta}(M_{Q}^{n-8} \times B^8)&=&(*,*,0,*,1/2^{4m-1}+3/2^{2m},1/2^{4m-2}+3/2^{2m})\\
\overrightarrow{\eta}(M_{Q}^{n-5})&=&(*,*,0,*,1/2^{4m+3}+3/2^{2m+2},1/2^{4m+2}+3/2^{2m+2})
\end{array}$$
Here $*$ is a term we aren't interested in, and we have divided by two already for the real representations. Thus in $\R/\Z$ we may consider the following matrix of vectors spanned by positive scalar curvature manifolds:
\begin{displaymath}
\left(\begin{array}{cccccc}
2^{-m-1} & 0 & 2^{-m-1} & * & * & * \\
0 &2^{-4m-3} & 2^{-4m-3} & 2^{-4m-3} & 2^{-4m-2} & 2^{-4m-3}\\
* & * & 0 & 2^{-2m-2} & 0 & 0 \\
* & * & 0 & * & 1/2^{4m+3}+3/2^{2m+2} & 1/2^{4m+2}+3/2^{2m+2} \\
* & * & 0 & * & 1/2^{4m-1}+3/2^{2m}  & 1/2^{4m-2}+3/2^{2m}\\
\end{array} \right)
\end{displaymath}
Note that since $\hat{D_8}$ and $\hat{C_8}$ both restrict to the same representation on $Q_8$, we can do a cancelation with the first three columns, namely, adding the first and the third and then subtracting the second gives us the column $(2^{-m},0,0,0,0)^T$, and so this, combined with the zeroes down the $\hat{Q_8}$ column immediately imply that we have spanned a subspace of order at least $2^m$ times what is spanned by $\RP^n$ and the quaternion lens spaces, which we already calculated as $2^{8+12m}$. Thus in total we have a subspace of order $2^{8+13m}$ as required.

Proceeding in exactly the same manner in $n=8m+7$ gives us a subspace of order $2^{11+13m}$ in $\R/\Z$. However here, the 'cancelation' we use above is slightly different, because when we restrict to the real representation of $C_8$ we have that $\eta(L^{8m+7}(1, \cdots,1))(\rho_{4}-\rho_0) \in \R/\Z$ in fact has the same order as $\eta(L^{8m+7}(1, \cdots,1))(2\rho_{4}-2\rho_0) \in \R/2\Z$, since $2\rho_4$ is in fact a quaternion representation. Thus we get an extra factor of $2$, giving the required $2^{12+13m}$ in total.

In dimensions $4m+1$ we use the Lens space bundles $L^{4m+1}$ described above, with fundamental group $C_8=<s>$, and calculate some eta invariants, again closely following the methods in \cite{2}. We note from \cite{2} that there is a surjective map from ${\L}_n(BC_8)$ to ${\L}_{n-4}(BC_8)$, where ${\L}_n(BC_8) \subset \R/\Z$ is the subspace spanned by the set of lens spaces $\eta(L^{n})(\rho)$, where $\rho \in R_0(C_l)$. Thus using naturality there is also a surjective map $\delta:{\L}_n(BSD_{16}) \rightarrow {\L}_{n-4}(BSD_{16})$.

Recall from the second section that in dimensions $8m+5$ we have $Ker(Ap)=[2^{m+1}]$, while in $8m+9$ we have $Ker(Ap)=[2^{m+1}]$ together with eta multiples. Since the order of the subgroups we need to realize using these lens space bundles is the same in $8m+5, 8m+9$, the above paragraph implies we need only realize a subgroup of order $2^{m+1}$ in dimensions $8m+5$. Further, again using \cite{2}, we need only check this claim in dimensions $5,13$, since then we would deduce that $\delta:{\L}_{13}(BSD_{16})=[4] \rightarrow {\L}_{9}(BSD_{16})=[2]$ is surjective with kernel $[2]$, and then by periodicity $\delta:{\L}_{8m+13}(BSD_{16}) \rightarrow {\L}_{8m+9}(BSD_{16})$ is also surjective with kernel at least $[2]$, so that inductively we would have realized a subgroup of order $2^{m+2}$ as required.

So we do some calculations in dimensions $5,13$. Note that $\chi_{\rho}$ restricts to $\rho_1 \oplus \rho_3$ on the cyclic subgroup $C_8=<s>$, where $\rho_1$ is the natural representation with $s \mapsto \omega=(1+i)/\sqrt{2}$, and $\rho_3$ with $s \mapsto \omega^3$. So we take:
$$\eta(i_*(L^5(1,1))(2\rho_0-\chi_{\rho})=\eta(L^5(1,1)(2\rho_0-\rho_1-\rho_3)$$
$$=\eta(L^5(1,1)(\rho_0-\rho_1)+\eta(L^5(1,1)(\rho_0-\rho_3)$$
and we calculate each summand directly using Lemma 3.4:
$$\eta(L^5(1,1)(\rho_0-\rho_1)=\frac{1}{8} \sum_{1 \neq \lambda \in C_8} \frac{\lambda(1+\lambda)(1-\lambda)}{(1-\lambda)^3}$$
$$=\frac{1}{8}\{\frac{\omega(1+\omega)}{(1-\omega)^2}+ \frac{i(1+i)}{(1-i)^2}+\frac{\omega^3(1+\omega^3)}{(1-\omega^3)^2}+0-\frac{\omega(1-\omega)}{(1+\omega)^2}-\frac{i(1-i)}{(1+i)^2}-
\frac{\omega^3(1-\omega^3)}{(1+\omega^3)^2}\}$$
$$=\frac{1}{8}\{\omega(\frac{(1+\omega)^3-(1-\omega)^3}{(1-i)^2})+\omega^3(\frac{(1+\omega^3)^3-(1-\omega^3)^3}{(1+i)^2})-\frac{(1+i)}{2}-\frac{(1-i)}{2}\}$$
$$=-3/4-1/8,$$
where we use $\omega^2=i, \omega^{4+j}=-\omega^j$, and then separate the $i,-i$ terms in the summand. We proceed analogously for $\rho_3$:
$$\eta(L^5(1,1)(\rho_0-\rho_3)=\frac{1}{8} \sum_{1 \neq \lambda \in C_8} \frac{\lambda(1+\lambda)(1-\lambda^3)}{(1-\lambda)^3}$$
$$=\frac{1}{8}\{\omega(\frac{(1+\omega)^4(1-\omega^3)-(1-\omega)^4(1+\omega^3)}{(1-i)^3})+\omega^3(\frac{(1+\omega^3)^4(1-\omega)-(1-\omega^3)^4(1+\omega)}{(1+i)^3})-\frac{2}{(1-i)^3}-\frac{2}{(1+i)^3}\}$$
$$=-3/4+1/8$$
So that adding up we get $3/2 \in \R/\Z$ of order 2, as required. We now make the analogous calculation in dimension $13$ for $L^{13}(1,1,1,1,1,1)$:
$$\eta(L^5(1,1,1,1,1,1)(\rho_0-\rho_1)=\frac{1}{8} \sum_{1 \neq \lambda \in C_8} \frac{\lambda^3(1+\lambda)(1-\lambda)}{(1-\lambda)^7}$$
$$=\frac{1}{8}\{\omega^3(\frac{(1+\omega)^7-(1-\omega)^7}{(1-i)^6})+\omega(\frac{(1+\omega^3)^7-(1-\omega^3)^7}{(1+i)^6})-\frac{i^3(1+i)}{(1-i)^6}+\frac{i(1-i)}{(1+i)^6}\}$$
$$=-17/8-1/32$$
While for $\rho_3$ we have:
$$\eta(L^5(1,1,1,1,1,1)(\rho_0-\rho_3)=\frac{1}{8} \sum_{1 \neq \lambda \in C_8} \frac{\lambda^3(1+\lambda)(1-\lambda^3)}{(1-\lambda)^7}$$
$$=\frac{1}{8}\{\omega^3(\frac{(1+\omega)^8(1-\omega^3)-(1-\omega)^8(1+\omega^3)}{(1-i)^7})+\omega(\frac{(1+\omega^3)^8(1-\omega)-(1-\omega^3)^8(1+\omega)}{(1+i)^7})+\frac{2}{(1-i)^7}+\frac{2}{(1+i)^7}\}$$
$$=-17/8+1/32$$
So that adding gives us $17/4$ of order $4 \in \R/\Z$, which completes the proof.
\end{proof}

\section{The missing class in $8k$ on the one column}
The result of the eta invariant calculations in the previous section, is that we can realize all of $Ker(Ap)$ in dimensions $4k+3$ and together with eta-multiples, in $4k+1$ also. However, in dimensions $8k$, there is a class on the one column of the local cohomology spectral sequence, which is not detected in periodic K-theory for which we need a separate geometric construction. This class is also detected in the ordinary $\Z_{2}$-homology of $BSD_{16}$, and its image is $\xi(yuP^{2k-1})$, the class dual to $yuP^{2k-1}$. We realize this class as follows:

The idea is to start with the exact sequence $C_{8} \rightarrow SD_{16} \rightarrow \Z_{2}$, lift it to obtain a sequence $C_{8} \rightarrow G \rightarrow \Z$, and take classifying spaces. Of course $S^{1} =B\Z$, and since $C_{8}$ acts on $n\C$ by multiplication, we have a map from a lens space $L$ into $BC_{8}$, and thus a fibre bundle $ L^{2n-1} \rightarrow M^{2n} \rightarrow S^{1}$, as shown. As a lens bundle this carries positive scalar curvature, and we claim that this gives the remaining class when $n \equiv 0 \mod 4$.
$$
\xymatrix{
L^{2n-1} \ar[r] \ar [d] & M^{2n} \ar[r] \ar [d]^{f} & S^{1} \ar[d]^{\simeq} \\
BC_{8} \ar[r] \ar [d] & BG \ar[r] \ar [d]^{g} & B\Z \ar[d] \\
BC_{8} \ar[r] & BSD_{16} \ar[r] & B\Z_{2} }
$$
\begin{prop}
When $n =4k$, the fundamental class of the manifold $M^{2n}$ constructed above has image $\xi(yuP^{2k-1}) \in H_*(BSD_{2^N})$, dual to $yuP^{2k-1} \in H^*(BSD_{2^N})$.
\end{prop}
\begin{proof}
Let $F=g \circ f$. We calculate the cohomology of the manifold $M$, along with $F^*:H^*(BSD_{16})\rightarrow H^*(M)$. We know $H^{*}(L^{2n-1}; \Z_{2})=\Z_{2}[X, \widetilde{\tau}]/(X^{n}, \widetilde{\tau}^{2})$ where $X$ is a dimension $2$ generator, and $\widetilde{\tau}$ dimension $1$, and $H^{*}(S^{1})=\Z_{2}[\sigma]/(\sigma^{2})$. By the universal coefficient theorem, all the cohomology groups except of course $H^{0}(M),H^{2n}(M)$ have rank $2$, and the Serre spectral sequence implies that there are two degree one and one degree two generators, one of the degree one generators being $\sigma$, with $\sigma^{2}=0$. Since $y=w_{1}(\chi_{2})$ has kernel $C_8$, restricting representations implies that $\sigma=F^{*}(y)$, and we can then define $\tau=F^{*}(x)$, and we will deduce  $H^{*}(M^{2n},\Z_{2})=\Z_{2}[\sigma, \tau, Z]/(\sigma^{2}, \sigma \tau + \tau^{2}, Z^{n})$, where $Z$ is a degree 2 class restricting to $X \in H^{*}(L^{4n-1})$.

Since $P$ restricts non-trivially on the cyclic group also, it follows that $F^*(P)=Z^2$ or $Z^2+Z\tau^2$. Either way $F^*(P^{2j})=Z^{4j}$, and recall that $P = c_{2}(\rho_1)$ reduced modulo $2$ and we can view $L^{2n-1}=S(n \sigma)/ C_{8}$. Notice that the universal cover $\widetilde{M}$ of $M=S(2k \rho_1) \underbrace{\times}_{G} \R$ is just $ S(2k \rho_1) \times \R$, and so the induced vector bundle $F^*(2k\rho_1)$ over $M$ is given by
 $$ 2k \rho_1 \underbrace{\times}_{G}(S(2k \rho_1) \times \R) \rightarrow  S(2k \rho_1) \underbrace{\times}_{G} \R $$
 This has a section via the diagonal map, which implies $Z^{n}=w_{2n}(F^*(2k\rho_1))$ must be zero. 

 Now $u^2=P(x^2+y^2) \in H^*(SD_{16})$ implies $F^*(u^2)=\tau^2 Z^2$, and since $xu=0$ we deduce $F^*(u)=Z(\tau + \sigma)$.
  We now claim this manifold $M$ is spin, and $F^*(P)=Z^2+Z\tau^2$. We have $Sq^1(u)=0$, so that $0=Sq^1(Z(\tau + \sigma))=Sq^1(Z)(\tau + \sigma)+Z\tau^2$ so that $Sq^1(Z)=Z\sigma$ or $Z(\tau + \sigma)$. Seeking a contradiction, assume the latter. Then $Sq^1(Z^3 \sigma)=Z^3\sigma \tau$ and $Sq^1(Z^3 \tau)=Z^3\tau^2$.

  We recall that in a smooth manifold N, $w_{k}(N)=\sum_{i+j=k} Sq^{i}(v_{j})$, where $v_{j} \in H^{j}(M^{n})$ is the unique class such that $v_{j}y=Sq^{j}(y)$, for every $y \in H^{n-j}(M)$. Thus we deduce from above that $w_1=v_1=\tau$. However, $M$ is a fibre bundle with orientable fibre, and so $w_1(M)$ must restrict to zero in $H^*(L^{4n-1})$. Since $\tau$ restricts to $\widetilde{\tau}$, we have a contradiction.

  So $Sq^1(Z)=Z\sigma$, implying $Sq^1(Z^3 \sigma)=Z^3\sigma^2=0$ and $Sq^1(Z^3 \tau)=Z^3\tau^2 + Z^3\sigma \tau=0$, so that $w_1(M)=v_1=0$ .
  Further $Sq^2(Z^3)=Z^4 +(Sq^1(Z))^2 =0$ and $Sq^2(Z^2 \tau^2)=\tau^2 Sq^2(Z^2)= \tau^2 (Sq^1(Z))^2=0$ so that $w_2(M)=v_2=0$ also.

  Further, $Sq^2(P)=u^2$ implies $Sq^2(F^*(P))=Z^2 \tau^2$, so that $F^*(P) \neq Z^2$.
 The top cohomology class of $M$ is $Z^{2n-1} \tau^{2}=Z^{n-1}\sigma \tau$, and observe that $F^*(yuP)=\sigma Z(\sigma +\tau)(Z^2+Z\tau^2)=Z^3\sigma \tau$. Any other class mapping to $Z^{2n-1}\sigma \tau$ must have a factor of $uP^{n-1}$ and as $xu=0$, dualising $F^{*}$ immediately gives $F_{*}(\xi(Z^{2n-1} \tau^{2})=\xi(yuP^{n-1})$.
 \end{proof}
 This construction is entirely analogous to the one used to obtain similar missing classes for dihedral groups. However, in that case, while the class in question is also not detected in periodic K-theory, it appears in the second local cohomology. This is a representation theoretic geometric construction, and the only immediate difference between the two cases is that the representation in question here is not real.

\section{ Ordinary cohomology and the two column}

The preceding two sections realize all of $Ker(Ap)$ that lies in the first local cohomology. It remains to consider the two column in the local cohomology spectral sequence. This is in fact the same for all semi-dihedral groups, and it is a $\Z_{2}$ vector space of dimension $k+1$ if $n=8k+4$,$k$ if $n=8k+j$ with $j \neq 4$ even, and $0$ if $n$ is odd. It suffices to detect these classes in ordinary $\Z_{2}$ homology. We claim that sufficiently many classes are realized by inclusion from a Klein four subgroup $V(2)$.

We refer here to \cite{8}, which gives explicit manifold constructions to calculate $H_{*}^{+}(BV(2))$ explicitly, together with the image under inclusion in $H_{*}^{+}(BD_8)$.

Note that the dihedral group of order $8$ includes into all the semi-dihedral groups.
We denote by $\omega \in D_{8}$ the rotation by $\pi/2$, $s,s',t,t'$ the reflections through the lines $y=0,x=0,y=x,y=-x$ respectively. The Klein four subgroup we consider is $V(2)=<s,s'>$ (note that $<t,t'>$ is another one, but it turns out not to be needed because viewed inside $SD_{16}$, it is conjugate to $<s,s'>$).

It is known that $H^{*}(BD_{8},\Z_{2})=\Z_{2}[\alpha,\beta,\delta]/(\alpha \beta +\beta^{2})$, where $\alpha=w_{1}(\hat{<\omega>}),\beta=w_{1}(\hat{<t,t'>}),\delta=w_{2}({\sigma})$ where $\sigma$ is a natural two-dimensional representation. Here $\hat{H}$ denotes the inflated representation, meaning the representation with kernel $H$,  and $\sigma$ is any natural two dimensional representation. Further we have $H^{*}(BV(2))=\Z_{2}[p,q]$, with $p=w_{1}(\hat{(\omega)^{2}}),q=w_{1}(\hat{s})$. Restricting representations then gives us the induced maps in cohomology. The map $H^{*}(BD_{8})\rightarrow H^{*}(BV(2))$ is determined by $ \alpha \rightarrow p, \beta \rightarrow 0, \delta \rightarrow q(p+q)$. So if $A=\Z_{2}[\alpha,\delta]$ we deduce $H^{*}(BV(2))=A \oplus qA$.
$$H^{*}(BD_{8},\Z_{2}) \rightarrow A \rightarrowtail H^{*}(BV(2))=A \oplus qA $$

 Following \cite{8}, the plan now is to dualise, and choose dual bases in homology and compute where the classes in $H_{*}^{+}(BV(2)$ ( which we know explicitly) map to.
 $$H_{*}(BD_{8},\Z_{2}) \leftarrow A^{v} \twoheadleftarrow H_{*}(BV(2))=A^{v} \oplus qA^{v} $$
 Denote the above inclusion by $f_*$. We now consider $\xi_{(a,b)} \in H_{n}(BV(2)$, dual to $p^a q^b$ in the monomial basis. Then using this decomposition $\xi_{(a,b)} \rightarrow \Sigma \xi(k)$ where $k \in A$ and $k \rightarrow p^{a}q^{b} + \cdots $, where $\xi$ denotes the dual basis to the basis $\{ \alpha^i \delta^j,\beta^i \delta^j\}$ of $H^{*}(BD_{8},\Z_{2})$, with $i,j \geq 0$. The following may be found in \cite{8}.
 \begin{prop}
The image under the inclusion $f_*$ of $H_{*}^{+}(BV(2))$ is spanned by the following classes:\\
In dimensions $4k+2$ it is spanned by $\xi(\alpha^{4i}\delta^{4j+3})$, with $0 \leq i,j \in Z$ and $4i+8j+6=4k+2$.\\
In dimensions $4k$ it is spanned by $\xi(\alpha^{4i+2}\delta^{4j+1})$, with $0 \leq i,j \in Z$ and $4i+2+8j+2=4k$.
\end{prop}
\begin{proof}
We start with $n=4k+2$. Then, using the dual basis to the monomial basis $H_{n}^{+}(BV(2))$ is $k-$ dimensional generated by $\xi_{(a,b)}$, the images of $\RP^{a} \times \RP^{b}$, with $a,b\equiv 3 \mod 4$.  We see that $\alpha^{4j}\delta^{a} \rightarrow p^{4k}q^{a}(p+q)^{a}=p^{4k+a}q^{a} + \cdots$ so that $\xi_{(a+4j,a)} \rightarrow \xi(\alpha^{4j} \delta^{a})$ (Note that $j=0$ is allowed here, with $a \equiv 3 \mod 4$). Since we start at $a=3$, adding these up gives us $\lfloor(k+1)/2 \rfloor$ classes in $H_{*}^{+}(BD_{8})$.

For $n=4k$ $H_{n}^{+}(BV(2)$ is $k+1-$ dimensional generated by $\xi_{(4k-1,1)},\xi_{(1,4k-1)}$, represented by $\RP^{4k-1} \times \RP^{1}$, and by all the $\xi_{(a,n-a)} +\xi_{(a-2,n-a+2)}$ which are the images of $\RP(2L_{0} \oplus (n-1-a)\varepsilon \rightarrow \RP^{a})$, with $ 5 \leq a \equiv 1 \mod 4$. Dualising again, we now check that the classes we obtain are exactly all those of the form $\xi(\alpha^{n-8i-2} \delta^{4i+1})$, with $0 \leq i \leq (n-2)/8$. Adding these gives $\lfloor(k+1)/2 \rfloor$ again.

Indeed, $\alpha^{n-2} \delta \rightarrow p^{n-1}q + \cdots $ so we get one class from the product $\RP^{n-1} \times \RP^{1}$. Recall that $[M_{(a,b)}]$ gives us the class $\xi_{(a,n-a)} +\xi_{(a-2,n-a+2)}$, so it will be mapped to $\Sigma \xi(k)$, where $k$ maps to a sum containing exactly one of $p^{a}q^{n-a}$ and $p^{a-2}q^{n-a+2}$. Note that when $n=12$, $\alpha^{2}\delta^{5} \rightarrow p^2q^{10}+p^3q^9+p^6q^6+p^7q^5$ so $[M_{(9,5)}]$ (and $[M_{(5,9)}]$) map to $\xi(\alpha^{2}\delta^{5})$. We now observe that if we have realized $\xi(\alpha^{i}\delta^{j}) \in H_{*}^{+}(BD_{8})$ this way, then we also realize $\xi(\alpha^{i+4k}\delta^{j}) \in H_{*}^{+}(BD_{8})$, because if $j=1$ it is just a product with a higher dimensional projective space, and otherwise the original class was realized by some $M_{(a,b)}$, so that $\xi(\alpha^{i+4k}\delta^{j})$ will be realized by one of $M_{(a+4,b)}$ or $M_{(a,b+4)}$. Similarly, if we have $\xi(\alpha^{2}\delta^{j}) \in H_{*}^{+}(BD_{8})$ then we also have $\xi(\alpha^{2}\delta^{j+4k})$. The $j=k=1$ case is the $n=12$ scenario done above, so just as before, if the original class was realized by some $M_{(a,b)}$, then $\xi(\alpha^{2}\delta^{j+4k})$ will be realized by one of $M_{(a+8,b)}$ or $M_{(a,b+8)}$.

So, to fully verify the proposition, we claim we can have no term of the form $\xi(\alpha^{i}\delta^{4j+3})$. Indeed, note that $n \equiv 0 \mod 4$ so $i \equiv 2 \mod 4$, and so it is clear that a product $\RP^{n-1} \times \RP^{1}$ can not hit this class. This leaves the possibility of some $M_{(a,b)}$, meaning that $\alpha^{i}\delta^{4j+3} \rightarrow p^{i}(q(p+q))^{4j+3}$ must map to $p^a,q^{b-2} + \cdots$ or $p^{a-2}q^b + \cdots$ but not $p^{a-2}q^{b}+p^a q^{b-2} + \cdots$. So equating binomial coefficients we must have $(4j+3)C_{a-i} \equiv 1 \mod 2$ and $(4j+3)C_{a-i-2} \equiv 0 \mod 2$, or vice versa.

However note that $a-i \equiv 3 \mod 4$, and for any $I,J$ it's clear that $(4J+3)C_{4I+3}$ and $(4J+3)C_{4I+1}$ differ by the factor $(4(J-I)+2)(4(J-I)+3)/((4I+2)(4I+3))$, which simplifies to odd numbers, so working modulo $2$, the two coefficients are the same.
\end{proof}

We thus have a sequence of inclusions $ko_{*}(BV(2)) \rightarrow ko_{*}(BD_{8}) \rightarrow ko_{*}(BSD_{16})$, and we have understood what classes are induced from the first inclusion. We know consider the second inclusion, starting with the cohomology of semi-dihedral groups and how they restrict. The following may be found in \cite{kijti}, \cite{evens}.\\
\begin{prop}\label{HSDN}
We have $$H^{*}(BSD_{N},\Z_{2})=\Z_{2}[x,y,u,P]/(xy+x^{2},xu,x^{3},u^{2}+(x^{2}+y^{2})P)$$
where $\mid x \mid=\mid y \mid=1,\mid u \mid=3,\mid P \mid=4 $. Further the restriction map $f^*:H^{*}(BSD_{N}) \rightarrow H^{*}(BD_{8})$ sends $x \rightarrow 0,y \rightarrow \alpha, u \rightarrow \alpha \delta, P \rightarrow \delta^2$.
\end{prop}
So we can immediately dualise again to check that we get sufficiently many classes in $H_{*}^{+}(BSD_{2^{N+2}};\Z_{2})$. We can now prove the following proposition, which implies the GLR conjecture for $SD_{16}$ immediately.
\begin{prop}
The image under the above inclusion spans all of the two-column in the local cohomology filtration for $ko_*(BSD_{2^{N+2}})$.
\end{prop}
\begin{proof}
We check that we have sufficiently many classes in ordinary homology.\\
We claim it suffices to check that if $i>0$, then $f_*(\xi(\alpha^i \delta^j))\neq 0$ and $f_*(\xi(\alpha^i \delta^j))=f_*(\xi(\alpha^{i'} \delta^{j'})) \Rightarrow i=i', j=j'$.\\
Indeed, we know that for $n=4k, 4k+2$, we have $\lfloor(k+1)/2 \rfloor$ classes in $H_*^+(BD_8)$ induced from $V(2)$. By Table $2.2$ we need to produce $K+1$ classes in $H_n^+(BSD_{2^{N+2}})$ for $n=8K+4$, and $K$ classes for $n=8K, 8K+2$ and $8K+6$.

Now if we have $8K+4=4k$ then $8K+8=4k+4$ which implies that $K+1=(k+1)/2=\lfloor(k+1)/2 \rfloor$. Similarly $8K=4k$ and $8K+2=4k+2$ both imply $K=k/2=\lfloor(k+1)/2 \rfloor$ since $k$ must be even. Finally if $8K+6=4k+2 \Rightarrow K=(k-1)/2=\lfloor(k+1)/2 \rfloor-1$, which is exactly what we need since in these dimensions we have an extra $\xi(\delta^{4K+3}) \in H_{*}^{+}(BD_{8},\Z_{2})$ which clearly maps to zero under $f_*$.

So it suffices to check that terms dual to $\alpha^i \delta^j$ with $i>0$ are mapped monomorphically by $f_*$. Recall that the classes $\xi(\alpha^i \delta^j) \in H_*^+(BD_8)$ always have $i$ even and $j$ odd, so that by Proposition $5.5.2$ we have that $f^*(y^{i-1}u P^{(j-1)/2})=\alpha^i \delta^j$. Thus $f_*(\xi(\alpha^i \delta^j)) \neq 0$, and we claim that $f_*(\xi(\alpha^i \delta^j))=\xi(y^{i-1}u P^{(j-1)/2})$.

This follows since if $f^*(y^a u^b P^c)=\alpha ^{a+b} \delta^{2c+b}$ and $f^*(y^{a'} u^{b'} P^{c'})=\alpha ^{a'+b'} \delta^{2c'+b'}$ both map to $\alpha^i \delta^j$ with $i>0$ even and $j \geq 1$ odd, then $a, a', b,b'$ must all be odd. Further $a+b=a'+b',2c+b=2c'+b'$, which implies that $a'=a-2j,b'=b+2j,c'=c-j$, for some $j \in \Z$, and $j \geq 0$ without loss of generality. However in $H^{*}(BSD_{2^{N+2}},\Z_{2})$, we have $u^2=x^2 P +y^2 P$ which implies $y u^3=yx^2uP +y^3 u P=y^3 u P$ since $xu=0$. Thus we can repeatedly apply this formula to deduce
$$y^{a'}u^{b'}P^{c'}=y^{a-2j}u^{b+2j}P^{c-j}=yu^3(y^{a-2j-1}u^{b+2j-3}P^{c-j})=y^3uP(y^{a-2j-1}u^{b+2j-3}P^{c-j})=$$
$$=y^{a-2j+2}u^{b+2j-2}P^{c-j+1}= \cdots =y^au^bP^c$$
Now since $f_*(\xi(\alpha^i \delta^j))= \sum \xi(k)$ where the sum is over all $k \in H^*(BSD_{2^{N+2}})$ that map under restriction to $\xi(\alpha^i \delta^j)$, we deduce that $f_*(\xi(\alpha^i \delta^j))=\xi(y^{i-1}u P^{(j-1)/2})$, and thus that $f_*(\xi(\alpha^i \delta^j))=f_*(\xi(\alpha^{i'} \delta^{j'})) \Rightarrow i=i', j=j'$, as required.
\end{proof}

\end{document}